\documentclass[a4paper]{article}

\pdfoutput=1 %preprint only

\usepackage[T1]{fontenc}
\usepackage[utf8]{inputenc}
\usepackage[english]{babel}
\usepackage{graphics}
\usepackage{graphicx}
\usepackage{epsfig}
\usepackage{amsmath}
\usepackage{amssymb}
\usepackage{amsthm}
\usepackage{url}
\usepackage{hyperref}
\usepackage{float}
\usepackage{pifont}
\usepackage[table]{xcolor}
\usepackage{lmodern}
\usepackage{xspace}
\usepackage{bbold}
\usepackage[capitalize,nameinlink]{cleveref}
\usepackage{subcaption}
\usepackage{mathtools}
\usepackage{adjustbox}

\definecolor{lightgray}{gray}{.85}

\usepackage{draftwatermark} %preprint only
\SetWatermarkText{Preprint} %preprint only
\SetWatermarkScale{1}       %preprint only

\newcommand{\R}{\mathbb{R}}

\newcommand{\T}{\intercal}
\newcommand{\Diff}{d}
\DeclareMathOperator{\e}{e}
\DeclareMathOperator{\tr}{tr}

\DeclareMathOperator{\vect}{vec}
\DeclareMathOperator{\sign}{sgn}

\title{Comparing (Empirical-Gramian-Based) \\ Model Order Reduction Algorithms}

\author{Christian Himpe\thanks{%
Computational Methods in Systems and Control Theory,
Max Planck Institute for \newline Dynamics of Complex Technical Systems,
Sandtorstr.~1, 39106 Magdeburg, Germany; \newline
ORCID:~0000-0003-2194-6754,
\url{himpe@mpi-magdeburg.mpg.de}
}}

\date{}

\begin{document}

\maketitle

\paragraph{Abstract} ~\\
In this work, the empirical-Gramian-based model reduction methods:
Empirical poor man's truncated balanced realization, empirical approximate balancing, empirical dominant subspaces, empirical balanced truncation,
and empirical balanced gains are compared in a non-parametric and two parametric variants,
via ten error measures: Approximate Lebesgue $L_0$, $L_1$, $L_2$, $L_\infty$, Hardy $H_2$, $H_\infty$, Hankel, Hilbert-Schmidt-Hankel,
modified induced primal, and modified induced dual norms, for variants of the thermal block model reduction benchmark. \linebreak
This comparison is conducted via a new meta-measure for model reducibility called MORscore.

\section{Introduction} % 0.5 pages
Model reduction research has made great strides in the past decades, spawning ever new methods and variants for specific requirements.
Yet, this plethora of algorithms is not (or only very sparsely) evaluated against each other on common benchmarks.
Such comparisons would enable a faster transfer of mathematically research to engineering and industrial applications. 

In the following, prototypically, a comparison of empirical-Gramian-based methods is demonstrated for a standard benchmark system in a manner,
which can be automated, for example to test various variants of a method determining the best suited for a problem.
In the scope of this work, model reduction for affine-parametric, generalized, linear time-invariant systems is considered:
\begin{align}\label{eq:linsys}
\begin{split}
 E \dot{x}(t) &= A(\theta)x(t) + Bu(t), \\
         y(t) &= Cx(t),
\end{split}
\end{align}
which consist of an ordinary differential equation in $x$, with a non-singular mass matrix $E \in \R^{N \times N}$,
an affinely decomposable parametric system matrix $A(\theta) = A_0 + \sum_{p=1}^P \theta_p A_p \in \R^{N \times N}$,
so that $E^{-1} A(\theta)$ is asymptotically stable for all parameters $\theta \in \R^P$,	
and an input matrix $B\in \R^{N \times M}$, as well as a linear output function defined by the output matrix $C \in \R^{Q \times N}$.

In the following some fundamentals of projection-based model reductions are assumed;
for a background on this topic the reader is referred to the seminal textbook~\cite{morAnt05}.

\pagebreak

\section{Empirical Gramians for Linear Systems} % 1-1.5 pages
System Gramians are system-theoretic operators encoding the input-output system properties of controllability and observability \cite{Kal60}.
Empirical Gramians \cite{morLalMG99} are generalizations of these system Gramians, which are based on quadrature,
and were introduced to apply linear, Gramian-based methods from linear system theory to nonlinear systems,
while incorporating nonlinear information and avoiding (explicit) linearization.
Since linear systems are a special case of nonlinear systems,
with, admittedly, a very simple ``nonlinearity'',
empirical Gramians can also be computed for linear systems.
Note that for linear systems, the empirical Gramians correspond to the classic system Gramians up to numerical error;
this is shown in \cite{morLalMG99,morHimO14}.
The quality of the empirical Gramians depends on simulated state and output trajectories for which the system is excited by perturbed input or initial state.
These perturbations are defined by scales ($c_m$ and $d_q$), which in this context are set to one, but in general should reflect the operating region of the system.
Following, we summarize the three fundamental empirical system Gramians in the case of linear systems.

\subsection{Empirical Controllability Gramian}
The controllability Gramian quantifies the ability to drive a linear system to a steady state in finite time via the input \cite{Kal63}.
For linear systems, the controllability Gramian matrix is defined as $W_C := \int_0^\infty \e^{E^{-1} At} E^{-1} B B^\T E^{-\T} \e^{A^\T E^{-\T} t} \Diff t$,
and classically computed as the (low-rank) solution to the Lyapunov equation $A W_C E^\T + E W_C A^\T = - BB^\T$.
Based on the definition of $W_C$, the empirical controllability Gramian is given by:
\begin{align*}
 \widehat{W}_C := \sum_{m=1}^M \int_0^\infty x^m(t) x^m(t)^\T \Diff t,
\end{align*}
with $x^m(t)$ being the solution of $E \dot{x}^m(t) = A x^m(t) + B (c_m e_m \delta(t))$,
suitable scales $c_m \in \R$, and the $m$-th canonical standard base vector $e_m \in \R^M$.

\subsection{Empirical Observability Gramian}
The observability Gramian matrix describes the ability to determine the state of linear system via the output in finite time \cite{Kal63}.
For linear systems, the observability Gramian matrix is defined as $W_O := \int_0^\infty \e^{A^\T E^{-\T} t} C^\T C \e^{E^{-1} At} \Diff t$,
and is classically computed as the (low-rank) solution to the Lyapunov equation $A^\T W_O E + E^\T W_O A = - C^\T C$.
Based on the definition of $W_O$, the (linear) empirical observability Gramian (via the dual system's controllability Gramian~\cite{morWilP02}) is given by:
\begin{align*}
 \widehat{W}_O := \sum_{q=1}^Q \int_0^\infty z^q(t) z^q(t)^\T \Diff t,
\end{align*}
with $z^q(t)$ being the solution of $E^\T \dot{z}^q(t) = A^\T z^q(t) + C^\T (d_q e_q \delta(t))$,
suitable scales $d_q \in \R$, and the $q$-th canonical standard base vector $e_q \in \R^N$.

\pagebreak

\subsection{Empirical Cross Gramian}
The cross Gramian matrix combines controllability and observability information and hence delineates the minimality of a linear system \cite{morFerN83}.
For \emph{square} linear systems (featuring the same number of inputs and outputs), the cross Gramian matrix $W_X$ is defined as $W_X := \int_0^\infty \e^{E^{-1} A t} E^{-1} B C \e^{E^{-1} At} \Diff t$,
and classically computed as the (low-rank) solution of the Sylvester equation \linebreak $A W_X E + E W_X A = -B C$.
Based on the definition of $W_X$, the (linear) empirical cross Gramian \cite{morBauBHetal17} is given by:
\begin{align*}
 \widehat{W}_X := \sum_{m=1}^M \int_0^\infty x^m(t) z^m(t)^\T \Diff t,
\end{align*}
with $x^m(t)$ being the solution of $E \dot{x}^m(t) = A x^m(t) + B (c_m e_m \delta(t))$,
$z^m(t)$ being the solution of $E^\T \dot{z}^m(t) = A^\T z^m(t) + C^\T (d_m e_m \delta(t))$,
suitable scales $c_m$, \linebreak $d_m \in \R$, and the $m$-th canonical standard base vector $e_m \in \R^M$.

For non-square systems, the non-symmetric cross Gramian $W_Z$,
the cross Gramian of the average system $(A, \bar{B} = \sum_{m=1}^M B_{*,m}, \bar{C} = \sum_{q=1}^Q C_{q,*}, E)$, is proposed in \cite{morHimO16}.
The linear empirical non-symmetric cross Gramian is given by:
\begin{align*}
 \widehat{W}_Z := \sum_{m=1}^M \sum_{q=1}^Q \int_0^\infty x^m(t) z^q(t)^\T \Diff t,
\end{align*}
with $x^m(t)$ being the solution of $E \dot{x}^m(t) = A x^m(t) + \bar{B} (c_m e_m \delta(t))$,
$z^m(t)$ being the solution of $E^\T \dot{z}^q(t) = A^\T z^q(t) + \bar{C}^\T (d_q \epsilon_q \delta(t))$,
suitable scales $c_m, d_q \in \R$, and the $m$-th, $q$-th canonical standard base vectors $e_m \in \R^M$, $\epsilon_q \in \R^Q$.

\subsection{Parametric Empirical Gramians}
Empirical Gramians may also be applied to parametric systems.
Here, the approach from \cite{morHimO15a} is utilized, which follows the general principle behind empirical Gramians: averaging over an operating region.
Hence, given a pre-selected sampling from parameter-space $\Theta_h$, an average (controllability, observability, cross, or non-symmetric cross) Gramian is computable \cite{morBauBHetal17}:
\begin{align*}
 \overline{W}_*(\Theta_h) := \sum_{\theta \in \Theta_h} W_*(\theta).
\end{align*}
For low-dimensional parameter-spaces, this could be some uniform grid in a region of interest;
for higher dimensional parameter-spaces, sparse grids can be utilized \cite{morBauB08b}.

Even though this averaging process can lead to annihilation,
it can be justified by the related accumulation process,
typically used, i.e., in (balanced) proper orthogonal decomposition (POD) model reduction \cite{morWilP02},
which (compresses and) concatenates trajectories before assembling a Gramian matrix.
So, given two discrete trajectory matrices $X_1$ and $X_2$, which are first concatenated and then a Gramian matrix is formed,
as for the abstract computation of a POD,
\begin{align*}
 \begin{bmatrix} X_1 & X_2 \end{bmatrix} \begin{bmatrix} X_1 & X_2 \end{bmatrix}^\T = X_1 X_1^\T + X_2 X_2^\T,
\end{align*}
this is mathematically (but not numerically due to annihilation) equivalent to the sum of the individual trajectory Gramians.

\pagebreak

\section{Empirical-Gramian-Based Model Reduction}\label{sec:methods} % 2 pages
Following, five empirical-Gramian-based model reduction methods are summarized,
of which either can be computed via the empirical controllability and observability Gramians $\{W_C, W_O\}$,
or via the empirical cross Gramian $W_X$ (empirical non-symmetric cross Gramian $W_Z$ for non-square systems).

The considered empirical-Gramian-based model reduction methods are exclusively projection-based approaches,
meaning from the empirical system Gramian matrices ``projection'' matrices are obtained --
a reducing projection $V$ and a reconstructing projection $U$, both of column-rank $n$:
\begin{align*}
 U \in \R^{N \times n}, \quad V \in \R^{n \times N},
\end{align*}
which appropriately applied to the system \eqref{eq:linsys} yield a reduced order system:
\begin{align*}
 (V E U) \dot{\tilde{x}}(t) &= \big((V A_0 U) + \sum_{p=1}^P \theta_p (V A_p U) \big) \tilde{x}(t) + (V B) u(t), \\
         \tilde{y}(t) &= (C U) \tilde{x}(t),
\end{align*}
or in a more compact form, as the reduced system matrices can be precomputed:
\begin{align*}
 \widetilde{E} \dot{\tilde{x}}(t) &= \widetilde{A}(\theta) \tilde{x}(t) + \widetilde{B} u(t), \\
         \tilde{y}(t) &= \widetilde{C} \tilde{x}(t).
\end{align*}
An orthogonal projection $U = V^\T$, $V U = I$ is called (Bubnov-)Galerkin projection,
a bi-orthogonal projection $U \neq V^\T$, $V U = I$ is called Petrov-Galerkin projection,
and a projection $U \neq V$, $V U \neq I$ is just called oblique projection.

In the following, only the features of the considered model reduction techniques are briefly summarized,
for a description and algorithm of these methods consult the referenced works in the respective subsections.
Note, that even though error bounds and error indicators are mentioned below for each method,
the purpose of this work is the heuristic comparison of methods against each other.

\subsection{Empirical Poor Man}
The Poor Man's Truncated Balanced Realization (PM) from \cite{morPhiS05} just utilizes either the (empirical) controllability Gramian,
or the (empirical) observability Gramian,
and uses the Gramian's dominant singular vectors as Galerkin projection.
Using the controllability Gramian in this fashion is equivalent to the proper orthogonal decomposition (POD),
using the observability Gramian is equivalent to the adjoint proper orthogonal decomposition~\cite{morBuiW05} (aPOD).

Being a Galerkin projection, this method is stability preserving in the reduced order model if the system is dissipative.
As an error indicator, typically the normalized sum of kept singular values is used as well as projection error of the data \cite{morOrSK12},
which quantifies the reduced model's preserved energy in relation to the full model.

\pagebreak

\subsection{Empirical Approximate Balancing}
Approximate balancing (AB) is a technique suggested in \cite[M3]{morRahVA14},
which uses the left and right singular vectors from a truncated SVD of the cross Gramian as oblique projection,
yet, without the bi-orthogonality of the Petrov-Galerkin projections, but orthogonality of the reducing and reconstructing projections with respect to themselves.
This method is based on the approximate balancing method from \cite{morSorA02}, but omits the eigenvector approximation. 
The counterpart variant based on controllability and observability Gramians is known as modified proper orthogonal decomposition \cite{morOrSK12},
which uses singular vectors from truncated SVDs of $W_C$ and $W_O$ similarly as oblique projection.
Even though, this method is claimed to be ``effective for non-normal systems'' (\cite[Sec.~III.D]{morOrSK12}),
for either method no error bounds or stability guarantees are available,
but as indicated in \cite[Fig.~8]{morOrSK12}, an error indicator can be derived based upon the projection error.
Due to the missing bi-orthogonality between the reducing and reconstructing projections,
it is paramount to apply the projections to the mass matrix if $E = I$.
Using empirical controllability, observability or cross Gramians yields the empirical approximate balancing method.

\subsection{Empirical Dominant Subspaces}
The dominant subspaces (DS) method, constructs a Galerkin projection by combining the dominant controllability and observability subspaces \cite{morPen06},
obtained from the respective (empirical) Gramians;
while the variant based on the (empirical) cross Gramian is introduced in \cite{morBenH19}.
The column-rank of the projection is then determined by the conjoined and orthogonalized singular vectors of the system Gramians,
weighted by their associated singular values.
As an orthogonal projection, DS is stability preserving for dissipative systems.
Furthermore, a Hardy-2 error bound exists for the controllability and observability Gramian-based DS \cite{morShiS05} (in two variants),
while a Lebesgue-$2$ error indicator is introduced in \cite{morBenH19} for the cross-Gramian-based DS.
To obtain and conjoin the system Gramians' singular vectors, various algorithms are available,
here, we use the truncated SVDs and rank-revealing SVDs for this task.

\subsection{Empirical Balanced Truncation}
Balanced truncation (BT) first transforms the system into a coordinate system in which controllability and observability are aligned,
via a Petrov-Galerkin projection, so the respective controllability and observability Gramians are diagonal and equal.
The diagonal entries, the Hankel singular values (HSVs), measure controllability and observability simultaneously,
hence the sub-system associated to the small HSVs is truncated.
This method from \cite{morMoo81} is the gold standard of system-theoretic model reduction methods,
due to, first, preserving stability in the reduced order model \cite{morPerS82},
and second, error bounds in the Hardy-$\infty$ norm \cite{morGlo84,morEnn84}, Hardy-$2$ norm \cite{morSorA02,morAnt05} and Lebesgue-$1$ norm \cite{morLamA92,morObiA01}.

To balance the Gramians $\{W_C, W_O\}$, the balanced POD ansatz \cite{morWilP02} is employed,
which corresponds to the square-root method \cite{morTomP87}, but using SVD-based square-roots of the Gramians.
Note that this does not lead to an exactly balanced system \cite[MR3]{morVar91b}.
For the $W_X$ ($W_Z$) balanced truncation variant, the method from \cite{morJiaQY19} is used, which in turn is based on \cite{morSafC88,morSafC89}.

\pagebreak

\subsection{Empirical Balanced Gains}
Balanced gains (BG) is a variant of balanced truncation,
of which the simplified variant from \cite{morDav86} is used here.
In balanced gains, the system is balanced as for balanced truncation,
but instead of the Hankel singular values, or the sum thereof,
an alternate measure is utilized, based on an observation on the $L_2$-norm of the impulse response (of symmetric systems):
\begin{align*}
 \|y\|_2^2 &= \tr(C W_C C^\T) = \tr(B^\T W_O B) = \tr(C W_X B) \\
            &= \sum_{k=1}^N \hat{c}_k^\T \hat{c}_k \sigma_k = \sum_{k=1}^N \hat{b}_k \hat{b}_k^\T \sigma_k = \sum_{k=1}^N |\hat{b}_k \hat{c}_k| \sigma_k,
\end{align*}
for the $k$-th row $\hat{b}_k$ of the balanced input matrix $\widehat{B}$,
and the $k$-th column $\hat{c}_k$ of the balanced output matrix $\widehat{C}$.
Hence, the sequence of base vectors is given by the magnitude of the quantity $d_k$, instead of the HSVs $\sigma_k$:
\begin{align*}
 d_k := \hat{c}_k^\T \hat{c}_k \sigma_k = \hat{b}_k \hat{b}_k^\T \sigma_k = |\hat{b}_k \hat{c}_k| \sigma_k.
\end{align*}
This means compared to balanced truncation, the same modes are used, but in a different order.
As the order of modes is not a requirement for stability preservation in the reduced order model,
it also holds for balanced gains \cite[Corollary~2]{morPerS82}.
Empirical balanced gains is then given by the (simplified) balanced gains approach using empirical Gramians.

\section{Approximate Norms}\label{sec:norms}
To comprehensively compare the reduced to the full order models,
four signal norms, four system norms, and two induced norms are applied.
For an elaborate discussion of these norms see \cite[Ch.~5,6]{BoyB91},\cite[Ch.~5]{morAnt05},\cite[Ch.~2]{Tos13}.
Due to numerical, efficiency or practical reasons, only approximate norms of the error system are considered.
Note, that the signal norms are computed from time-domain trajectories,
and the system (and modified induced) norms are approximated by transformations of empirical Gramians,
instead of frequency domain sampling.

\subsection{Signal Norms}
The signal norms are based on time-domain evaluations of the system output $y$ and the reduced system's output $\tilde{y}$,
and are given as the Lebesgue norms of the output error $\|y - \tilde{y}\|$.
Practically, vector norms of vectorized discrete output trajectories $y_h$, $\tilde{y}_h$ ($Q$ outputs $\times$ $K$ time steps data matrices) are computed.

\subsubsection{Approximate $L_0$-``Norm''}
The $L_0$ signal ``norm'' describes the sparsity of a \emph{discrete-time} signal \cite{SchEA11},
and is approximated, based on \cite{arash13}, for an error signal by:
\begin{align*}
 \|y_h - \tilde{y}_h\|_{L_0} = \sum_{k=0}^K \sum_{q=1}^Q \big|\sign\big(y_{h,q}(k)-\tilde{y}_{h,q}(k)\big)\big| \approx \sqrt[n]{\prod_{\ell=1}^{QK} |\vect(y_h - \tilde{y}_h)_\ell}|.
\end{align*}
Technically, this is not a norm, due to the lack of absolute scalability, but for the intended purpose this function can be treated as a norm.

\pagebreak

\subsubsection{Approximate Lebesgue $L_1$-Norm}
The Lebesgue $L_1$-norm of a signal quantifies the action or consumption of a process and its definition and approximation for an output error signal are given by:
\begin{align*}
 \|y - \tilde{y}\|_{L_1} = \int_0^\infty \|y(t)-\tilde{y}(t)\|_1 \Diff t \approx \Delta t \, \| \vect(y_h - \tilde{y}_h) \|_1;
\end{align*}
in terms of the model reduction error it can also be seen as the area under the error signal.

\subsubsection{Approximate Lebesgue $L_2$-Norm}
The Lebesgue $L_2$-norm of a signal measures its energy.
Its definition and approximation for an output error signal are given by:
\begin{align*}
 \|y - \tilde{y}\|_{L_2} = \sqrt{\int_0^\infty \|y(t)-\tilde{y}(t)\|_2^2 \Diff t} \approx \sqrt{\Delta t} \, \| \vect(y_h - \tilde{y}_h) \|_2,
\end{align*}
which can be interpreted as the energy loss in the reduced order model.
As all methods tested in this work are energy-based, this norm is the canonical error measure.

\subsubsection{Approximate Lebesgue $L_\infty$-Norm}
The Lebesgue $L_\infty$-norm of a signal determines its peak,
with definition and approximation of the error signal given by:
\begin{align*}
 \|y - \tilde{y}\|_{L_\infty} = \sup_t \|y(t) - \tilde{y}(t)\|_\infty \approx \| \vect(y_h - \tilde{y}_h) \|_\infty,
\end{align*}
which yields the maximum error between the signals.

\subsection{System Norms}
The system norms characterize frequency-domain errors of the reduced system's output 
$G_r(\omega) := C_r(E_r \omega-A_r)^{-1}B_r$ compared to the system output \linebreak $G(\omega) := C(E\omega-A)^{-1}B$,
for frequencies $\omega \in \mathbb{C}$, $\operatorname{Re}(\omega) < 0$, and are either Hardy-norms and/or Schatten-norms of the Hankel operator $H$.
These four norms were selected based on \cite[Sec.~2.2.7]{morSch96}.

\subsubsection{Approximate Hardy $H_2$-Norm}
The Hardy $H_2$-norm can be interpreted as the root-mean-square of the frequency response to white noise,
the $L_2$-norm of the impulse response (thus also known as impulse response norm),
the maximum output amplitude for finite input, or average gain.
To approximate the $H_2$-norm, the truncated balanced part of the output operator and controllability Gramian are utilized \cite[Remark~3.3]{morSorA02}:
\begin{align*}
 \|G - G_r\|_{H_2} = \sqrt{\int \tr((G(\imath \omega)-G_r(\imath \omega)) (G(\imath \omega)-G_r(\imath \omega))^*) \Diff \omega} \approx \sqrt{\widehat{\bar{C}}_2 W_{Z,22} \widehat{\bar{B}}_2}.
\end{align*}

\subsubsection{Approximate Hardy $H_\infty$-Norm}
The Hardy $H_\infty$-norm describes the worst-case frequency domain error,
which relates, via Parseval's equation, to the maximum $L_2$-gain,
and thus to the time-domain $L_2$ error.
Based on \cite[Corollary~9.3]{morGlo84}, the $H_\infty$ error can be approximated by the balanced truncation error bound,
which in turn is approximated by the principal discarded Hankel singular value \cite[Ch.~2.4]{morHim17}:
\begin{align*}
 \|G - G_r\|_{H_\infty} = \sup(\sigma_1(G(\imath \omega)-G_r(\imath \omega))) \approx 2 \!\! \sum_{k=n+1}^N \sigma_k(H) \approx 2 (N - n) \sigma_{n+1}(H),
\end{align*}
and is related to the nuclear norm (Schatten-1 norm) of the Hankel operator.
Alternatively, the $H_\infty$-norm could be approximated by the trace of the non-symmetric cross Gramian $\|G - G_r\|_{H_\infty} \approx -\frac{1}{2} \tr(W_{Z,22}) = -\bar{C}_2 A_{22}^{-1} \bar{B}_2$ \cite{morLiuST98a}.

\subsubsection{Approximate Hilbert-Schmidt-Hankel-Norm}
The Hilbert-Schmidt-Hankel norm corresponds to the operator norm \linebreak (Schatten-2 norm) of the Hankel operator, and as for the $H_\infty$-norm,
is approximated using only the principal discarded Hankel singular value:
\begin{align*}
 \|G - G_r\|_{HSH} = \sqrt{\sum_{k=n+1}^N \sigma_k^2(H)} \approx \sqrt{(N-n) \sigma_{n+1}^2(H)}
\end{align*}
Scaled by a factor of $\pi$, the square-root of this norm yields the enclosed area of the Nyquist plot \cite{Han92}.

\subsubsection{Approximate Hankel-Norm}
The Hankel norm is given by the principal discarded singular value of the Hankel operator,
which corresponds to the Schatten-$\infty$ norm of the Hankel operator:
\begin{align*}
 \|G - G_r\|_{Ha} = \sigma_{n+1}(H).
\end{align*}
This norm is the lower bound for the model reduction error as by the Adamjan-Arov-Krein theorem \cite{morGlo87,morGloP87}.

\subsection{Modified Induced Norms}
If the Hankel operator is used in its classic form, it maps from and to a function space of squarely integrable functions,
and the (previous) Hankel norm is its induced norm.
If one modifies the Hankel operator to allow for a function space of just continuous functions as domain or range,
the induced norms change as follows \cite{morWil85}.
Note, that for single-input-single-output systems, the following norms coincide with the Hardy-2 norm.

\subsubsection{Induced Primal Norm}
Modifying the Hankel operator to the expanded domain of continuous functions,
the induced norm becomes the square-root of the input-observability Gramian's spectral radius:
\begin{align*}
 \|H - H_r \|_{\mathcal{H}_C} = \sqrt{\lambda_{max}(B_{22}^\T W_{O,22} B_{22})}.
\end{align*}

\subsubsection{Induced Dual Norm}
Modifying the Hankel operator to the expanded range of continuous functions,
is equivalent to expanding the dual system's Hankel operator's domain,
thus the induced norm becomes the square-root of the output-controllability Gramian's spectral radius:
\begin{align*}
 \|H - H_r\|_{\mathcal{H}_O} = \sqrt{\lambda_{max}(C_{22} W_{C,22} C_{22}^\T)}.
\end{align*}

\subsection{Parametric Norms}
To obtain an error quantification for parametric systems, the previous norms are extended with respect to the considered system's parameter-space.
Given a (state-space) error norm $\|\cdot\|_X$, the associated parametric state-space error norm is given by the composition with a parameter-space norm $\|\cdot\|_Y$.
In \cite{morBauBBetal11} (see also \cite{morBenGW15}), this composite state-parameter norms are defined via a norm as a mapping $\|\cdot\|_{X \otimes Y} : M \times \Theta \to \R_+$,
with the Cartesian product of output, response or operator domain $M$ and parameter domain $\Theta$ respectively.
To approximate these parametric norms, a sampling of the parameter-space $\Theta_h \subset \Theta$ is drawn,
and given this finite, discrete parameter sample $\Theta_h$ an approximate norm is computed.
We follow \cite{morGruHKetal13}, in evaluating the parametric $X \otimes L_1$, $X \otimes L_2$, and $X \otimes L_\infty$ norms:
\begin{align*}
 \|y(\theta) - \tilde{y}(\theta)\|_{X \otimes L_1} &= \;\;\;\;\int_\Theta \|y(\theta) - \tilde{y}(\theta)\|_X \Diff \theta &&\approx \;\;\;\,\sum_{\theta \in \Theta_h} \|y(\theta) - \tilde{y}(\theta)\|_X, \\
 \|y(\theta) - \tilde{y}(\theta)\|_{X \otimes L_2} &= \sqrt{\int_\Theta \|y(\theta) - \tilde{y}(\theta)\|_X^2 \Diff \theta} &&\approx \sqrt{\sum_{\theta \in \Theta_h} \|y(\theta) - \tilde{y}(\theta)\|_X^2}, \\
 \|y(\theta) - \tilde{y}(\theta)\|_{X \otimes L_\infty} &= \;\,\max_{\theta \in \Theta} \|y(\theta) - \tilde{y}(\theta)\|_X &&\approx \;\;\;\,\max_{\theta \in \Theta_h} \|y(\theta) - \tilde{y}(\theta)\|_X,
\end{align*}
for $X$ being any of the signal, system or induced norms.
To estimate the quality of a parametric reduced order model fairly, it is a basic requirement to have disjoint training and test parameter sets.
Typically, this is implicitly ensured by a (sparse) grid parameter sampling for the training and randomly drawn test parameters from a suitable distribution.

\section{MORscore}
The comparison of model reduction errors for varying reduced orders, see for example \cref{fig:err},
is a useful vehicle to evaluate the performance of model reduction techniques for a specific system in a certain norm.
Yet, there are multiple relevant features in these error graphs characterizing the associated model order reduction algorithm,
such as: lowest attained error or fastest error decay.
Now, a one-by-one comparison for multiple methods, in various norms is too tedious for potentially many systems.
A similar problem arises in comparing optimization codes,
which is managed by so-called \emph{relative minimization profiles} (RMP) \cite[Sec.~5]{CurMO17}.
These RMPs standardize such comparisons in various measures, such as best computed objective, and inspired the following scoring.
To make many-way model reduction comparisons feasible, a scalar score is introduced next,
summarizing a method's features in a specific norm based on the error graph.

\vskip1ex
\textbf{Definition}~(MORscore) \\ 
\textit{
Given an error graph $(n,\varepsilon(n)) \in \mathbb{N}_{>0} \times (0,1]$,
relating a reduced order $n$ to a relative output error of a model reduction method $M$ for a system $\Sigma$ in norm $\|\cdot\|$,
the normalized error graph $(\varphi_n,\varphi_{\varepsilon(n)})$ is determined by the maximum reduced order $n_{\text{max}} \in \mathbb{N}_{>0}$,
and machine precision $\epsilon_{\text{mach}} \in (0,1] \subset \R$ via mappings:}
\begin{align*}
 \varphi_n: \mathbb{N}_{>0} &\to [0,1], \;\; n \mapsto \frac{n}{n_{\text{max}}}, \\
 \varphi_\varepsilon: (0,1] &\to [0,1], \;\; \varepsilon \mapsto \frac{\log_{10}(\varepsilon)}{\lfloor\log_{10}(\epsilon_{\text{mach}})\rfloor},
\end{align*}
\textit{and the \textbf{MORscore} $\mu$ is defined as the area under this normalized error graph,}
\begin{align*}
 \mu_{(n_{\text{max}},\epsilon_{\text{mach}})}(M,\Sigma,\|\cdot\|) &:= \operatorname{area}(\varphi_n,\varphi_\varepsilon).
\end{align*}

By $\varphi_n$ the discrete reduced orders $1, 2 \dots n_{\text{max}}$ are mapped to the real interval $[0,1]$ by normalization.
And by $\varphi_\varepsilon$ the \emph{relative} model reduction error $\varepsilon$ is mapped to the real interval $[0,1]$,
by normalizing the $10$-base logarithm of the error by the $10$-base logarithm of the maximum accuracy $\epsilon_{\text{mach}}$ of the utilized number system;
i.e.\ double precision floating point numbers have an accuracy of approximately $\epsilon_{\text{mach}}(\text{dp}) \approx 10^{-16}$,
so $\lfloor\log_{10}(\epsilon_{\text{mach}}(\text{dp}))\rfloor = -16$.
Practically, the area is computed via the trapezoid rule\footnote{\url{https://www.mathworks.com/help/matlab/ref/trapz.html}}.
Note, that the maximum tested reduced order $n_{\text{max}}$ should be (far) below the original model order,
since the error decay flattens at some reduced order.
Hence, given a system of large order, and two model reduction methods, both yielding their minimal error reduced models at low orders,
a MORscore up to the full order would show only little difference.
Selecting the largest reduced order which attains the minimal error as $n_{\text{max}}$, the MORscore is a lot more meaningful.

Altogether, the MORscore is specified by the normalization, and describes the model reduction performance of a method for a system in a norm by single number,
as typical for (desktop) computer performance benchmarks.
A larger MORscore $\mu \in (0,1)$ means better model reduction performance, since the more area covered, the faster and lower the error decay.

As opposed to the $\beta$-RMPs \cite[Def.~5.2]{CurMO17}, no computational budget is prescribed here,  
nonetheless, the MORscore could be extended in this manner by limited computational time or even a prescribed $n_{\text{max}}$.

\section{Benchmark Comparison}
For a thorough comparison, the presented empirical-Gramian-based model reduction methods are tested in ten (approximate) norms for different configurations of a benchmark system.
In coordination with the model reduction software projects:
pyMOR \cite{morMilRS16}, MORLAB \cite{morBenW19b}, M.E.S.S \cite{SaaKB19-mmess-2.0}, a thermal block benchmark is tested.
A summary of the components for this comparison is given below.

\paragraph{Methods} ~\\
Each of the five methods summarized in \cref{sec:methods}, can be computed via the empirical controllability and observability Gramians $\{W_C,W_O\}$,
or the empirical (non-symmetric) linear cross Gramian $W_Z$.
Hence overall, ten empirical-Gramian-based model reduction techniques are compared: 
\begin{itemize}
 \item Empirical Poor Man (PM), via $W_C$ or $W_O$,
 \item Empirical Approximate Balancing (AB), via $\{W_C, W_O\}$ or $W_Z$,
 \item Empirical Dominant Subspaces (DS), via $\{W_C, W_O\}$ or $W_Z$,
 \item Empirical Balanced Truncation (BT), via $\{W_C, W_O\}$ or $W_Z$,
 \item Empirical Balanced Gains (BG), via $\{W_C, W_O\}$ or $W_Z$.
\end{itemize}

\paragraph{Parameterization} ~\\
In \cref{sec:benchmark}, a parametric benchmark with a four dimensional parameter-space is tested.
The benchmark is compared in three configurations:
\begin{itemize}
 \item Non-Parametric (parameters treated as constants),
 \item Single Parameter (parameters treated as single parameter),
 \item Multiple Parameters (parameters treated separately).
\end{itemize}

\paragraph{Measures} ~\\
The model reduction methods are compared via their MORscore for varying reduced orders in the following norms from \cref{sec:norms}:
\begin{itemize}
 \item Approximate Lebesgue $L_0$-``norm'',
 \item Approximate Lebesgue $L_1$-norm,
 \item Approximate Lebesgue $L_2$-norm,
 \item Approximate Lebesgue $L_\infty$-norm,
 \item Approximate Hardy $H_2$-norm,
 \item Approximate Hardy $H_\infty$-norm,
 \item Approximate Hilbert-Schmidt-Hankel-norm,
 \item Approximate Hankel-norm,
 \item Approximate modified induced primal norm,
 \item Approximate modified induced dual norm,
\end{itemize}
as well as the number of unstable ROMs up to the maximum order (denoted by the symbol $\mathcal{L}$).
Lyapunov stability is assessed via the real-part of the largest real eigenvalue of the pencil $(\widetilde{E},\widetilde{A}(\theta))$.
In the parametric case, these counts are averaged, similar to the considered norms, in an $L_1$, $L_2$ and $L_\infty$ sense over the sampled parameters.

\subsection{\texttt{emgr} -- EMpirical GRamian Framework} % 0.5 pages
All tested methods are based on empirical system Gramian matrices.
To compute these empirical Gramians for the subsequent numerical experiments,
the empirical Gramian framework \texttt{emgr} \cite{morHim18b} is adopted,
which has a unified interface \cite{morHimO13} for the empirical controllability, observability and (linear) cross Gramians.
Furthermore, the convergence of the empirical Gramians to the classic algebraic Gramians for linear systems is shown in \cite{morHim17}.
Practically, the current version \texttt{emgr} 5.7 \cite{morHim19a} is used.

\pagebreak

\subsection{Thermal Block Benchmark}\label{sec:benchmark} % 0.5 pages
For the comparison of the empirical-Gramian-based model order reduction methods, a recurring benchmark example (due to the well reducible diffusion process),
modeling the heat equation on the unit-square \cite[Thermal~Block]{morwiki} is utilized.

This thermal block benchmark system models dynamic heating of a two-\linebreak dimensional,
square domain $\Omega = (0,1) \times (0,1)$ with four enclosed circular regions $\omega_{i=1\dots 4}$ of equal radius,
one per quadrant, and each of individual parametric heat conductivity (diffusivity) $\kappa(x)$.
The left boundary of the domain $\partial \Omega_1 := \{0\} \times (0,1)$ is the inflow, realized by a Neumann boundary condition,
the top and bottom boundaries $\partial \Omega_2 := (0,1) \times \{0\}$, $\partial \Omega_4 := (0,1) \times \{1\}$ are insulated, via zero Neumann conditions,
while the right boundary $\partial \Omega_3 := \{1\} \times (0,1)$ prescribes Dirichlet-zero boundary conditions.
Lastly, the four quantities of interests $\mathcal{Y}_i$ are the average temperature of each circle $\omega_i$.
The overall partial differential equation (PDE) system is thus given by:
\begin{align*}
 \partial_t u(x,t) &= -\kappa(x) \Delta_x u(x,t), \quad &&x \in \Omega, \\
            \partial_x u(x,t) &= F(x,t), &&x \in \partial \Omega_1, \\
            \partial_x u(x,t) &= 0, &&x \in \partial \Omega_2 \cup \partial \Omega_4, \\
            u(x,t) &= 0, &&x \in \partial \Omega_3, \\
    \mathcal{Y}_i(t) &= \smallint_{\omega_i} u(x,t) \Diff x, \\
            \kappa(x) &= \begin{cases} \theta_i & x \in \omega_i, \;i=1\dots 4, \\ \theta_0 & \text{otherwise.} \end{cases}
\end{align*}
This PDE is discretized in space using the finite element method (FEM),
via the FEniCs software package \cite{fenics15}, yielding an ordinary differential equation system of the form \eqref{eq:linsys}.
The resulting linear input-output system has one input and four outputs, while the state-space has dimension $7488$,
and the parameter-space is four-dimensional, with $\theta_{i=1\dots 4} \in [1,10] \subset \R$ as in \cite{morBalK16},
while the background diffusivity constant is set to $\theta_0 = 1$.
For more a detailed description of this benchmark, and the software stack used for its creation, see also Chapter~\emph{(TBD)}.

\subsection{Numerical Results} % 4 pages

In the following, three variants of the thermal block benchmark are tested:
\begin{enumerate}
 \item No parameter: $\frac{1}{5}\theta_1 = \frac{2}{5}\theta_2 = \frac{3}{5}\theta_3 = \frac{4}{5}\theta_4 \equiv \sqrt{10}$,
 \item One parameter: $\frac{1}{5}\theta_1 = \frac{2}{5}\theta_2 = \frac{3}{5}\theta_3 = \frac{4}{5}\theta_4 \in [1,10]$,
 \item Four parameters: $\theta \in [1,10]^4$.
\end{enumerate}
For the parametric variants, the ($3 \cdot \dim(\theta)$) training samples of the parameter-space are taken from a logarithmically uniform grid,
whereas (ten) test samples are drawn randomly from a logarithmically uniform distribution over the parameter range.
The empirical Gramians are build from trajectories excited by impulses, while the ROMs are tested by random input.
The decompositions for the empirical-Gramian-based model reduction methods are approximated up to rank one-hundred.
Practically, the following numerical results are conducted using MATLAB 2019b on an Intel(R) Core(TM) i3-7130U CPU @ 2.70GHz with 8GB RAM.

\pagebreak

\subsubsection{Fixed Parameter}\label{sec:fixed}

\begin{figure}
 \includegraphics[width=\textwidth]{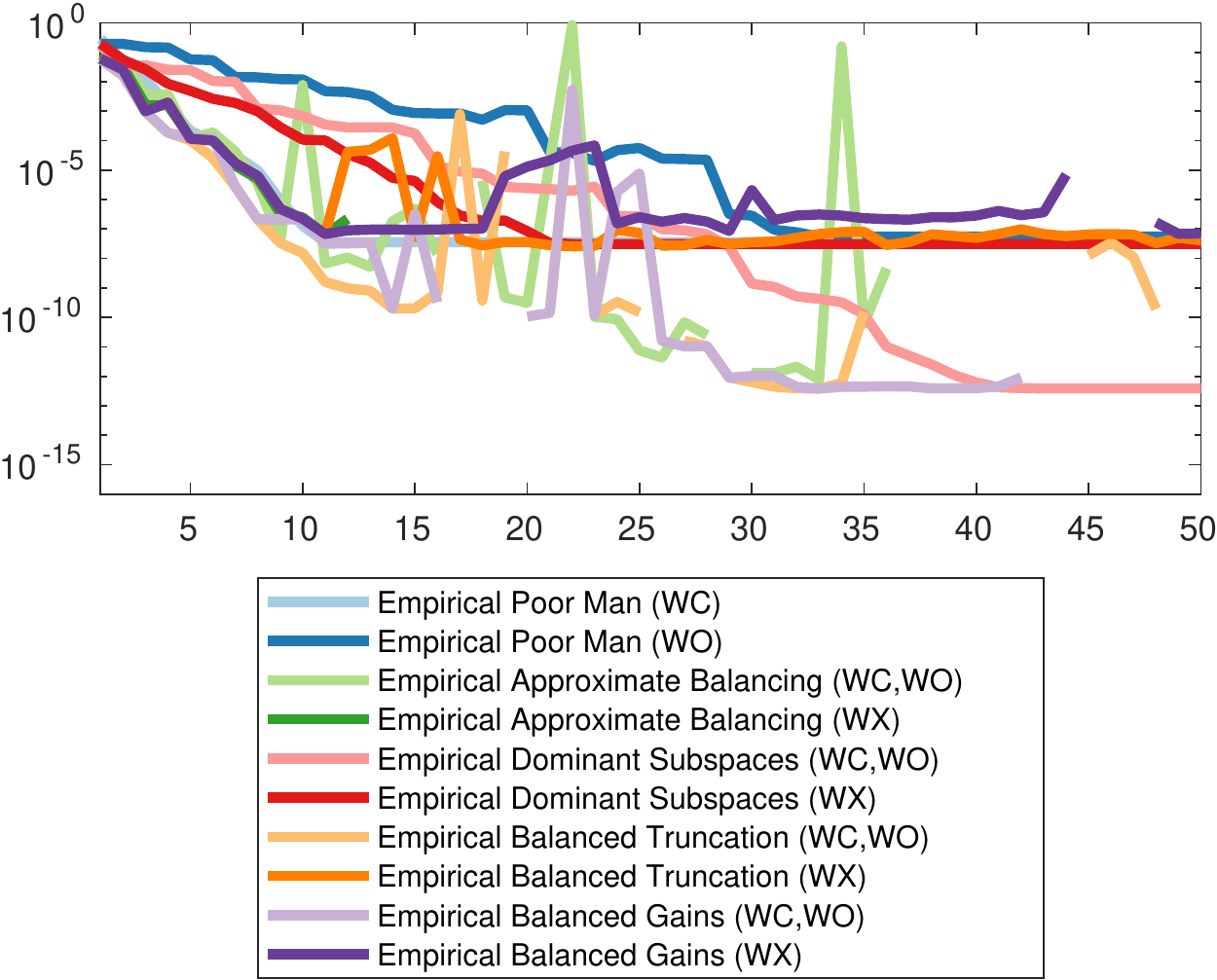}
 \caption{Relative error of reduced order models in the $L_2$-norm compared to the full order model for varying reduced orders.}
 \label{fig:err}
\end{figure}

In the first set of numerical experiments, the thermal block benchmark is tested with a single fixed parameter.
Exemplary in Figure~\ref{fig:err}, the model reduction error in the approximate $L_2$-norm for the ten considered methods are compared for reduced models of orders one to fifty.
This figure illustrates how complex a visualization already in a single norm is. 
The proposed MORscores are listed in Table~\ref{tab:para0}, which is similarly not directly decipherable by a human observer,
yet, algorithmically it can be processed.
In the approximate signal norms the maximum MORscores are achieved by the DS($W_C$,$W_O$),
closely followed by BG($W_C$,$W_O$).
Notably the BT variants used are not in lead,
which in this case is related to many unstable reduced order models,
originating in the low-rank approximation of the Gramians and using an SVD-based square-root method for balancing,
nullifying the stability-preservation of the original balanced truncation method.
While the Galerkin methods do not produce unstable ROMs, all Petrov-Galerkin methods produce at least $21$ unstable ROMs.
The $H_2$-norm is lead by the PM($W_C$) method, whereas the $H_\infty$, $Ha$ and $HSH$ norms are headed by DS($W_C$,$W_O$),
closely followed by PM($W_C$).
Finally, in modified induced norms $\mathcal{H}_C$ and $\mathcal{H}_O$,
PM($W_C$) and PM($W_O$) perform best respectively.
Overall for this benchmark, the methods using $W_C$ and/or $W_O$ outperformed methods using $W_Z$,
likely due to the non-square system, which requires additional averaging in the non-symmetric cross Gramian.

\begin{table}
\rowcolors{1}{white}{lightgray}
\begin{adjustbox}{center}
\begin{tabular}{l|ccccccccccc}
 & $L_0$ & $L_1$ & $L_2$ & $L_\infty$ & $H_2$ & $H_\infty$ & $H\!S\!H$ & $Ha$ & $\mathcal{H}_C$ & $\mathcal{H}_O$ & $\mathcal{L}$ \\
 \hline
 PM($W_C$)       & 0.42 & 0.42 & 0.41 & 0.39 & 0.63 & 0.49 & 0.51 & 0.52 & 0.54 & 0.06 & 0 \\
 PM($W_O$)       & 0.29 & 0.29 & 0.29 & 0.28 & 0.10 & 0.38 & 0.38 & 0.38 & 0.10 & 0.45 & 0 \\
 AB($W_C$,$W_O$) & 0.33 & 0.33 & 0.32 & 0.30 & 0.46 & 0.03 & 0.04 & 0.04 & 0.44 & 0.39 & 37 \\
 AB($W_X$)       & 0.08 & 0.08 & 0.08 & 0.08 & 0.35 & 0.02 & 0.02 & 0.02 & 0.35 & 0.04 & 38 \\
 DS($W_C$,$W_O$) & 0.45 & 0.45 & 0.44 & 0.43 & 0.32 & 0.51 & 0.52 & 0.52 & 0.29 & 0.25 & 0 \\
 DS($W_X$)       & 0.39 & 0.38 & 0.38 & 0.36 & 0.34 & 0.39 & 0.39 & 0.39 & 0.34 & 0.08 & 0 \\
 BT($W_C$,$W_O$) & 0.38 & 0.38 & 0.37 & 0.35 & 0.43 & 0.36 & 0.36 & 0.36 & 0.43 & 0.18 & 25 \\
 BT($W_X$)       & 0.41 & 0.40 & 0.39 & 0.38 & 0.28 & 0.30 & 0.30 & 0.30 & 0.28 & 0.08 & 21 \\
 BG($W_C$,$W_O$) & 0.43 & 0.43 & 0.42 & 0.41 & 0.42 & 0.35 & 0.35 & 0.35 & 0.42 & 0.17 & 25 \\
 BG($W_X$)       & 0.36 & 0.35 & 0.34 & 0.32 & 0.28 & 0.30 & 0.30 & 0.30 & 0.28 & 0.08 & 37
\end{tabular}
\end{adjustbox}
\caption{MORscores($50$,$\epsilon_{\text{mach}}(DP)$) for the non-parametric benchmark .}
\label{tab:para0}
\end{table}

\pagebreak

\subsubsection{Single Parameter}\label{sec:single}
The MORscores for the single parameter benchmark are given in Table~\ref{tab:para11} ($L_1$), Table~\ref{tab:para12} ($L_2$) and Table~\ref{tab:para18} ($L_\infty$).
Generally, all methods perform worse compared to the non-parametric benchmark,
since the averaging of empirical Gramians over parameter samples decreases specific accuracy while increasing general applicability.
The signal norms are lead by BT($W_C$,$W_O$) and directly followed by BG($W_C$,$W_O$), PM($W_C$), DS($W_C$,$W_O$), and DS($W_X$).
In the $H_2$ and $\mathcal{H}_C$ norms, the methods BT($W_C$,$W_O$), PM($W_C$), and AB($W_C$,$W_O$) are in the lead,
while in the system norms $H_\infty$, $HSH$, $Ha$ the PM($W_C$) heads the MORscores.
The $\mathcal{H}_O$ norm is topped by PM($W_O$) and AB($W_C$,$W_O$) methods.
Balanced gains (BG) seem to work well for this benchmark, while approximate balancing (AB) perform worst overall.
As for the non-parametric benchmark, the Galerkin methods consistently produce stable ROMs, and the Petrov-Galerkin methods tend to assemble unstable ROMs.

\subsubsection{Multiple Parameters}\label{sec:multi}
The MORscores for the multiple parameter benchmark are given in Table~\ref{tab:para41} ($L_1$), Table~\ref{tab:para42} ($L_2$) and Table~\ref{tab:para48} ($L_\infty$),
and correspond overall to the single parameter setting, yet, with again slightly lower scores.
Curiously, balanced gains performance drops more than balanced truncation.

\subsubsection{MORscore Discussion}
Summarizing, the presented MORscore tables can improve heuristic comparisons of model reduction methods.
An automated evaluation could include filtering extreme values per norm, as demonstrated in the previous evaluations,
or means per methods across norms.
Specifically for the comparison of the empirical-Gramian-based model reduction methods on the thermal block benchmark,
the arithmetic means of MORscores across norms yields the PM($W_C$) and DS($W_C$,$W_O$) methods
as top scoring for the non-parametric benchmark,
and the PM($W_C$) = POD for the parametric benchmark variants, as in \cite{morBauBHetal17}.

Beyond this sample comparison, the proposed MORscore could find application in model reduction software development signaling regressions,
or defining highscore boards of competing methods for benchmark problems.

\pagebreak

\section{Conclusion} % 1 page
This work should be considered an exemplary quantitative comparison using MORscores, and by no means exhaustive comparison.
Specifically, other relevant (empirical) Gramian-based methods not tested here are (empirical) singular perturbation approximation \cite{morFerN82a},
and (empirical) Hankel norm approximation \cite{morFerDF10},
yet both methods are not purely projection based,
but require a numerically potentially expensive post-processing of a balanced realization.
Also, the empirical Gramians have various variants \cite{morHim18b} that could be tested,
as well as different balancing algorithms \cite{morVar91b}.
Nevertheless, this work can serve as a template for benchmarking model reduction methods by their \textbf{MORscore}.

\section*{Code Availability Section}
The source code of the presented numerical examples can be obtained from:
\begin{center}
\url{http://runmycode.org/companion/view/3760}
\end{center}
and is authored by: \textsc{Christian Himpe}.

\section*{Acknowledgement}
Supported by the German Federal Ministry for Economic Affairs and
Energy (BMWi), in the joint project: ``MathEnergy -- Mathematical
Key Technologies for Evolving Energy Grids'', sub-project:
Model Order Reduction (Grant number: 0324019\textbf{B}).

\bibliographystyle{plainurl}
\bibliography{mor,csc,software,extra}

\pagebreak

\section*{Appendix}

\begin{table}[h!]\small
\rowcolors{1}{white}{lightgray}
\begin{adjustbox}{center}
\begin{tabular}{l|ccccccccccc}
 & $L_0$ & $L_1$ & $L_2$ & $L_\infty$ & $H_2$ & $H_\infty$ & $H\!S\!H$ & $Ha$ & $\mathcal{H}_C$ & $\mathcal{H}_O$ & $\mathcal{L}$ \\
 \hline
 PM($W_C$)       & 0.26 & 0.25 & 0.25 & 0.23 & 0.37 & 0.42 & 0.44 & 0.44 & 0.37 & 0.07 & 0 \\
 PM($W_O$)       & 0.18 & 0.18 & 0.18 & 0.17 & 0.10 & 0.23 & 0.24 & 0.24 & 0.10 & 0.18 & 0 \\
 AB($W_C$,$W_O$) & 0.15 & 0.15 & 0.14 & 0.14 & 0.35 & 0.03 & 0.04 & 0.04 & 0.36 & 0.18 & 37.5 \\
 AB($W_X$)       & 0.06 & 0.06 & 0.06 & 0.06 & 0.24 & 0.02 & 0.02 & 0.02 & 0.23 & 0.05 & 38.1 \\
 DS($W_C$,$W_O$) & 0.24 & 0.23 & 0.23 & 0.22 & 0.19 & 0.30 & 0.31 & 0.32 & 0.19 & 0.15 & 0 \\
 DS($W_X$)       & 0.24 & 0.23 & 0.23 & 0.22 & 0.24 & 0.29 & 0.29 & 0.30 & 0.24 & 0.07 & 0 \\
 BT($W_C$,$W_O$) & 0.25 & 0.25 & 0.24 & 0.24 & 0.38 & 0.28 & 0.28 & 0.28 & 0.36 & 0.14 & 14.8 \\
 BT($W_X$)       & 0.18 & 0.18 & 0.18 & 0.17 & 0.20 & 0.19 & 0.19 & 0.19 & 0.20 & 0.10 & 33.2 \\
 BG($W_C$,$W_O$) & 0.26 & 0.26 & 0.26 & 0.25 & 0.34 & 0.23 & 0.23 & 0.23 & 0.33 & 0.12 & 18.5 \\
 BG($W_X$)       & 0.12 & 0.12 & 0.12 & 0.11 & 0.19 & 0.18 & 0.18 & 0.18 & 0.19 & 0.08 & 34.2
\end{tabular}
\end{adjustbox}
\caption{MORscores($50$,$\epsilon_{\text{mach}}(DP)$) for the single parameter benchmark ($L_1$).}
\label{tab:para11}
\end{table}

\begin{table}[h!]\small
\rowcolors{1}{white}{lightgray}
\begin{adjustbox}{center}
\begin{tabular}{l|ccccccccccc}
 & $L_0$ & $L_1$ & $L_2$ & $L_\infty$ & $H_2$ & $H_\infty$ & $H\!S\!H$ & $Ha$ & $\mathcal{H}_C$ & $\mathcal{H}_O$ & $\mathcal{L}$ \\
 \hline
 PM($W_C$)       & 0.22 & 0.22 & 0.22 & 0.20 & 0.34 & 0.39 & 0.40 & 0.41 & 0.34 & 0.04 & 0 \\
 PM($W_O$)       & 0.15 & 0.15 & 0.15 & 0.14 & 0.07 & 0.20 & 0.21 & 0.21 & 0.07 & 0.15 & 0 \\
 AB($W_C$,$W_O$) & 0.11 & 0.11 & 0.10 & 0.10 & 0.32 & 0.00 & 0.01 & 0.01 & 0.33 & 0.15 & 118.66 \\
 AB($W_X$)       & 0.03 & 0.03 & 0.03 & 0.02 & 0.21 & 0.00 & 0.00 & 0.00 & 0.20 & 0.02 & 120.56 \\
 DS($W_C$,$W_O$) & 0.20 & 0.20 & 0.20 & 0.19 & 0.16 & 0.27 & 0.28 & 0.29 & 0.16 & 0.12 & 0 \\
 DS($W_X$)       & 0.20 & 0.20 & 0.20 & 0.19 & 0.21 & 0.26 & 0.26 & 0.27 & 0.21 & 0.04 & 0 \\
 BT($W_C$,$W_O$) & 0.21 & 0.21 & 0.21 & 0.20 & 0.35 & 0.25 & 0.25 & 0.25 & 0.33 & 0.10 & 47.03 \\
 BT($W_X$)       & 0.14 & 0.14 & 0.14 & 0.13 & 0.17 & 0.16 & 0.16 & 0.16 & 0.17 & 0.07 & 105.00 \\
 BG($W_C$,$W_O$) & 0.23 & 0.22 & 0.22 & 0.21 & 0.30 & 0.20 & 0.20 & 0.20 & 0.30 & 0.09 & 58.52 \\
 BG($W_X$)       & 0.09 & 0.09 & 0.08 & 0.08 & 0.16 & 0.15 & 0.15 & 0.15 & 0.16 & 0.05 & 108.16
\end{tabular}
\end{adjustbox}
\caption{MORscores($50$,$\epsilon_{\text{mach}}(DP)$) for the single parameter benchmark ($L_2$).}
\label{tab:para12}
\end{table}

\begin{table}[h!]\small
\rowcolors{1}{white}{lightgray}
\begin{adjustbox}{center}
\begin{tabular}{l|ccccccccccc}
 & $L_0$ & $L_1$ & $L_2$ & $L_\infty$ & $H_2$ & $H_\infty$ & $H\!S\!H$ & $Ha$ & $\mathcal{H}_C$ & $\mathcal{H}_O$ & $\mathcal{L}$ \\
 \hline
 PM($W_C$)       & 0.24 & 0.23 & 0.23 & 0.21 & 0.37 & 0.42 & 0.44 & 0.44 & 0.37 & 0.07 & 0 \\
 PM($W_O$)       & 0.17 & 0.17 & 0.17 & 0.16 & 0.10 & 0.23 & 0.24 & 0.24 & 0.10 & 0.18 & 0 \\
 AB($W_C$,$W_O$) & 0.12 & 0.12 & 0.12 & 0.11 & 0.35 & 0.03 & 0.04 & 0.04 & 0.36 & 0.18 & 40 \\
 AB($W_X$)       & 0.05 & 0.05 & 0.05 & 0.05 & 0.24 & 0.02 & 0.02 & 0.02 & 0.23 & 0.05 & 41 \\
 DS($W_C$,$W_O$) & 0.22 & 0.22 & 0.21 & 0.20 & 0.19 & 0.30 & 0.31 & 0.32 & 0.19 & 0.15 & 0 \\
 DS($W_X$)       & 0.22 & 0.22 & 0.22 & 0.21 & 0.24 & 0.29 & 0.29 & 0.30 & 0.24 & 0.07 & 0 \\
 BT($W_C$,$W_O$) & 0.23 & 0.23 & 0.22 & 0.21 & 0.38 & 0.28 & 0.28 & 0.28 & 0.36 & 0.14 & 17 \\
 BT($W_X$)       & 0.16 & 0.16 & 0.16 & 0.14 & 0.20 & 0.19 & 0.19 & 0.19 & 0.20 & 0.10 & 34 \\
 BG($W_C$,$W_O$) & 0.24 & 0.24 & 0.24 & 0.23 & 0.34 & 0.23 & 0.23 & 0.23 & 0.33 & 0.12 & 19 \\
 BG($W_X$)       & 0.10 & 0.10 & 0.10 & 0.09 & 0.19 & 0.18 & 0.18 & 0.18 & 0.19 & 0.08 & 35
\end{tabular}
\end{adjustbox}
\caption{MORscores($50$,$\epsilon_{\text{mach}}(DP)$) for the single parameter benchmark ($L_\infty$).}
\label{tab:para18}
\end{table}

\pagebreak

\begin{table}[h!]\small
\rowcolors{1}{white}{lightgray}
\begin{adjustbox}{center}
\begin{tabular}{l|ccccccccccc}
 & $L_0$ & $L_1$ & $L_2$ & $L_\infty$ & $H_2$ & $H_\infty$ & $H\!S\!H$ & $Ha$ & $\mathcal{H}_C$ & $\mathcal{H}_O$ & $\mathcal{L}$ \\
 \hline
 PM($W_C$)       & 0.24 & 0.23 & 0.23 & 0.22 & 0.30 & 0.33 & 0.34 & 0.35 & 0.29 & 0.08 & 0 \\
 PM($W_O$)       & 0.18 & 0.17 & 0.17 & 0.16 & 0.10 & 0.24 & 0.24 & 0.24 & 0.10 & 0.18 & 0 \\
 AB($W_C$,$W_O$) & 0.12 & 0.12 & 0.11 & 0.11 & 0.31 & 0.03 & 0.04 & 0.04 & 0.29 & 0.18 & 43.4 \\
 AB($W_X$)       & 0.09 & 0.08 & 0.08 & 0.08 & 0.18 & 0.02 & 0.02 & 0.02 & 0.18 & 0.07 & 33.0 \\
 DS($W_C$,$W_O$) & 0.21 & 0.21 & 0.20 & 0.19 & 0.20 & 0.30 & 0.32 & 0.33 & 0.20 & 0.16 & 0 \\
 DS($W_X$)       & 0.19 & 0.19 & 0.19 & 0.18 & 0.20 & 0.24 & 0.25 & 0.25 & 0.21 & 0.09 & 0 \\
 BT($W_C$,$W_O$) & 0.24 & 0.24 & 0.24 & 0.23 & 0.30 & 0.22 & 0.22 & 0.22 & 0.30 & 0.20 & 5.1 \\
 BT($W_X$)       & 0.08 & 0.08 & 0.08 & 0.07 & 0.15 & 0.14 & 0.14 & 0.14 & 0.15 & 0.11 & 29.8 \\
 BG($W_C$,$W_O$) & 0.20 & 0.20 & 0.20 & 0.19 & 0.27 & 0.19 & 0.19 & 0.19 & 0.27 & 0.18 & 7.6 \\
 BG($W_X$)       & 0.05 & 0.05 & 0.05 & 0.05 & 0.13 & 0.12 & 0.12 & 0.12 & 0.13 & 0.11 & 36.7
\end{tabular}
\end{adjustbox}
\caption{MORscores($50$,$\epsilon_{\text{mach}}(DP)$) for the multi parameter benchmark ($L_1$).}
\label{tab:para41}
\end{table}

\begin{table}[h!]\small
\rowcolors{1}{white}{lightgray}
\begin{adjustbox}{center}
\begin{tabular}{l|ccccccccccc}
 & $L_0$ & $L_1$ & $L_2$ & $L_\infty$ & $H_2$ & $H_\infty$ & $H\!S\!H$ & $Ha$ & $\mathcal{H}_C$ & $\mathcal{H}_O$ & $\mathcal{L}$ \\
 \hline
 PM($W_C$)       & 0.20 & 0.20 & 0.19 & 0.19 & 0.27 & 0.30 & 0.31 & 0.32 & 0.26 & 0.05 & 0 \\
 PM($W_O$)       & 0.14 & 0.14 & 0.14 & 0.13 & 0.07 & 0.21 & 0.21 & 0.21 & 0.07 & 0.15 & 0 \\
 AB($W_C$,$W_O$) & 0.08 & 0.08 & 0.07 & 0.07 & 0.28 & 0.00 & 0.01 & 0.01 & 0.25 & 0.15 & 137.59 \\
 AB($W_X$)       & 0.04 & 0.04 & 0.04 & 0.04 & 0.15 & 0.00 & 0.00 & 0.00 & 0.15 & 0.04 & 104.58 \\
 DS($W_C$,$W_O$) & 0.18 & 0.17 & 0.17 & 0.16 & 0.17 & 0.27 & 0.29 & 0.30 & 0.17 & 0.13 & 0 \\
 DS($W_X$)       & 0.16 & 0.15 & 0.15 & 0.15 & 0.17 & 0.21 & 0.22 & 0.22 & 0.18 & 0.05 & 0 \\
 BT($W_C$,$W_O$) & 0.20 & 0.20 & 0.20 & 0.19 & 0.27 & 0.19 & 0.19 & 0.19 & 0.27 & 0.17 & 16.76 \\
 BT($W_X$)       & 0.04 & 0.04 & 0.04 & 0.03 & 0.12 & 0.11 & 0.11 & 0.11 & 0.12 & 0.08 & 94.24 \\
 BG($W_C$,$W_O$) & 0.17 & 0.16 & 0.16 & 0.15 & 0.24 & 0.16 & 0.16 & 0.16 & 0.24 & 0.15 & 25.18 \\
 BG($W_X$)       & 0.01 & 0.01 & 0.01 & 0.01 & 0.10 & 0.08 & 0.09 & 0.09 & 0.10 & 0.08 & 116.19
\end{tabular}
\end{adjustbox}
\caption{MORscores($50$,$\epsilon_{\text{mach}}(DP)$) for the multi parameter benchmark ($L2$).}
\label{tab:para42}
\end{table}

\begin{table}[h!]\small
\rowcolors{1}{white}{lightgray}
\begin{adjustbox}{center}
\begin{tabular}{l|ccccccccccc}
 & $L_0$ & $L_1$ & $L_2$ & $L_\infty$ & $H_2$ & $H_\infty$ & $H\!S\!H$ & $Ha$ & $\mathcal{H}_C$ & $\mathcal{H}_O$ & $\mathcal{L}$ \\
 \hline
 PM($W_C$)       & 0.21 & 0.21 & 0.21 & 0.20 & 0.30 & 0.33 & 0.34 & 0.35 & 0.29 & 0.08 & 0 \\
 PM($W_O$)       & 0.16 & 0.16 & 0.16 & 0.15 & 0.10 & 0.24 & 0.24 & 0.24 & 0.10 & 0.18 & 0 \\
 AB($W_C$,$W_O$) & 0.09 & 0.09 & 0.09 & 0.09 & 0.31 & 0.03 & 0.04 & 0.04 & 0.29 & 0.18 & 47 \\
 AB($W_X$)       & 0.06 & 0.06 & 0.06 & 0.05 & 0.18 & 0.02 & 0.02 & 0.02 & 0.18 & 0.07 & 36 \\
 DS($W_C$,$W_O$) & 0.19 & 0.19 & 0.19 & 0.18 & 0.20 & 0.30 & 0.32 & 0.33 & 0.20 & 0.16 & 0 \\
 DS($W_X$)       & 0.17 & 0.17 & 0.17 & 0.16 & 0.20 & 0.24 & 0.25 & 0.25 & 0.21 & 0.09 & 0 \\
 BT($W_C$,$W_O$) & 0.22 & 0.21 & 0.21 & 0.20 & 0.30 & 0.22 & 0.22 & 0.22 & 0.30 & 0.20 & 9 \\
 BT($W_X$)       & 0.05 & 0.05 & 0.05 & 0.05 & 0.15 & 0.14 & 0.14 & 0.14 & 0.15 & 0.11 & 30 \\
 BG($W_C$,$W_O$) & 0.18 & 0.18 & 0.17 & 0.17 & 0.27 & 0.19 & 0.19 & 0.19 & 0.27 & 0.18 & 11 \\
 BG($W_X$)       & 0.02 & 0.02 & 0.02 & 0.02 & 0.13 & 0.12 & 0.12 & 0.12 & 0.13 & 0.11 & 39
\end{tabular}
\end{adjustbox}
\caption{MORscores($50$,$\epsilon_{\text{mach}}(DP)$) for the multi parameter benchmark ($L_\infty$).}
\label{tab:para48}
\end{table}

\end{document}